\documentclass[a4paper,11pt]{amsart}
\usepackage{amsmath,amsthm}
\allowdisplaybreaks
\newtheorem{theorem}{Theorem}[section]
\newtheorem{cor}[theorem]{Corollary}
\newtheorem{lemma}[theorem]{Lemma}
\newtheorem{prop}[theorem]{Proposition}
\theoremstyle{remark}
\newtheorem{remark}[theorem]{Remark}
\theoremstyle{definition}
\newtheorem{definition}[theorem]{Definition}
\numberwithin{equation}{section}
\DeclareMathOperator{\id}{id}
\DeclareMathOperator{\tr}{tr}
\DeclareMathOperator{\slice}{sl}
\DeclareMathOperator{\clsp}{\overline{span}}
\newcommand{\bs}{{\mathfrak B}}
\newcommand{\multi}{\Ww}
\newcommand{\abs}[1]{\lvert#1\rvert}
\newcommand{\norm}[1]{\lVert#1\rVert}
\newcommand{\cstar}{$C^*$\ndash}
\newcommand{\Star}{${}^*$\ndash}
\newcommand{\rankone}[2]{#1 \otimes\overline{#2}}
\newcommand{\set}[1]{\{#1\}}
\newcommand{\setspace}[1]{\{\,#1\,\}}
\newcommand{\ip}[2]{\langle #1, #2 \rangle}
\newcommand{\Ra}{\Rightarrow}
\newcommand{\ndash}{\nobreakdash-}
\newcommand{\bh}{\Bb(\Hh)}
\newcommand{\Fe}{F_\Ee}
\newcommand{\Bfe}{\Bb(\Fe)}
\newcommand{\Pb}{\Pp_\beta}
\newcommand{\Wb}{\Ww_\beta}
\newcommand{\To}{{\Tt\Oo}_n}
\newcommand{\F}{{\Ff}_n}
\newcommand{\be}{\mathbf e}
\newcommand{\field}[1]{\mathbb{#1}}
\newcommand{\CC}{\field{C}}
\newcommand{\NN}{\field{N}}
\newcommand{\TT}{\field{T}}
\newcommand{\ZZ}{\field{Z}}
\newcommand{\Bb}{{\mathcal B}}

\newcommand{\Ee}{{\mathcal E}}
\newcommand{\Ff}{{\mathcal F}}

\newcommand{\Hh}{{\mathcal H}}
\newcommand{\Jj}{{\mathcal J}}
\newcommand{\Kk}{{\mathcal K}}
\newcommand{\Oo}{{\mathcal O}}
\newcommand{\Pp}{{\mathcal P}}

\newcommand{\Tt}{{\mathcal T}}
\newcommand{\Uu}{{\mathcal U}}
\newcommand{\Ww}{{\mathcal W}}
\begin{document}
%
% Top Matter
%
\title{States of Toeplitz-Cuntz algebras}
\author{Neal J. Fowler}
\address{Department of Mathematics  \\
      University of Newcastle\\  Callaghan, NSW  2308 \\ AUSTRALIA}
\email{neal@maths.newcastle.edu.au}
\begin{abstract}
We characterize the state space of a Toeplitz-Cuntz algebra $\To$
in terms of positive operator matrices $\Omega$ on Fock space which
satisfy $\slice\Omega \le \Omega$, where $\slice\Omega$ is the operator matrix
obtained from $\Omega$ by taking the trace in the last variable.
Essential states correspond to those matrices $\Omega$ which are
slice-invariant.  As an application we show that a pure essential product
state of the fixed-point algebra for the action of the gauge group
has precisely a circle of pure extensions to $\To$.
\end{abstract}
\maketitle
%
% Document body
%
\section*{Introduction}
Let $\To$ be the unital \cstar algebra which is universal
for collections of $n$ isometries with mutually orthogonal ranges;
we call $\To$ a {\em Toeplitz-Cuntz} algebra.
Since their introduction by Cuntz (\cite{cun}),
these algebras have been profitably used in the study of
normal \Star endomorphisms of $\bh$, the algebra of bounded
operators on a Hilbert space $\Hh$.  The main idea is 
as follows.  Let $\set{v_k: 1 \le k \le n}$ be the distinguished generating
isometries in $\To$.  Every \Star representation
$\pi$ of $\To$ on $\Hh$ on gives rise to
an endomorphism $\alpha$ of $\bh$ via
\[
\alpha(A) = \sum_{k=1}^n \pi(v_k)A\pi(v_k)^*,\qquad A\in\bh,
\]
and every endomorphism is of this form for some $n$ and $\pi$;
see \cite{arv,laca1,bjp}.

Arveson has generalized these ideas to the continuous case through
the use of product systems.
Every representation $\phi$ of a continuous product system $E$ on Hilbert space
gives rise to a semigroup $\alpha = \setspace{\alpha_t:t>0}$
of endomorphisms of $\bh$,
and $\phi$ is said to be {\em essential\/} if each $\alpha_t$ is unital;
such semigroups are called {\em $E_0$-semigroups,\/} and are the primary
objects of study in Arveson's series \cite{arv,arv2,arv3,arv4}.

One of Arveson's key results is that every product system $E$ has an essential
representation.  To prove this he associates with $E$
a universal \cstar algebra $C^*(E)$ whose representations are in bijective
correspondence with representations of $E$, characterizes the state
space of $C^*(E)$, and then uses this characterization
to show that there are always certain states, called {\em essential\/} states,
whose GNS representations give rise to essential representations of $E$.

In this paper we develop a discrete version of Arveson's method
in which $E$ is a product system over the positive integers $\NN$.
The algebras which arise as $C^*(E)$ are precisely the Toeplitz-Cuntz algebras:
up to isomorphism there is a unique product system $E^n$
over $\NN$ for each $n\in\set{1,2,\dots,\infty}$,
and $C^*(E^n) \cong \To$; see \cite{fowrae}.
We write $\Tt\Oo_\infty$ for the Cuntz algebra $\Oo_\infty$,
a notation which underlies an important advantage of our methods:
they apply for both finite and infinite $n$, so that one
does not have to study $\Oo_\infty$ as a special case.
While our methods are motivated by those of Arveson,
our exposition avoids any explicit use of product systems:
since the \cstar algebra being analyzed is a familiar one,
we can use it as a starting point rather than the product system.

Our main result is Theorem~\ref{theorem:TOn isomorphism},
which characterizes the state space of $\To$ in terms of a class
of positive linear functionals on the \Star algebra $\bs$
of operators on Fock space which have ``bounded support'';
these functionals are the analogues of Arveson's decreasing locally normal
weights (\cite{arv4}).
In Theorem~\ref{theorem:TOn isomorphism2} we give a reformulation of
this result in terms of the so-called {\em density matrices\/} associated
with these functionals; these are certain infinite operator matrices
on full Fock space over an $n$\ndash dimensional Hilbert space.
Roughly speaking, positive linear functionals on $\To$ correspond
to positive matrices $\Omega$ of trace-class operators with the property
that
\[
\slice\Omega \le \Omega,
\]
where $\slice\Omega$, the {\em slice\/} of $\Omega$,
is the operator matrix obtained from $\Omega$ by
``taking the trace in the last variable.''

Viewing $\To$ as the universal \cstar algebra of a product system,
a state $\rho$ of $\To$ is essential if its associated GNS representation
$\pi$ satisfies
\[
\sum \pi(v_k)\pi(v_k)^* = I.
\]
In Proposition~\ref{prop:sing-ess} we use the results of \S1 to
give some alternate characterizations of both essentiality
and the complementary notion of singularity;
perhaps the most useful aspect of this theorem is the characterization of
essentiality in terms of invariance under the map $\beta^*$
defined by
\[
\beta^*\rho(v_{i_1}\dotsm v_{i_k} v_{j_l}^*\dotsm v_{j_1}^*)
 =\sum_{m=1}^n \rho(v_{i_1}\dotsm v_{i_k} v_m v_m^* v_{j_l}^*\dotsm v_{j_1}^*).
\]
When $n$ is finite, essential states of $\To$
are precisely those which factor through the canonical homomorphism
of $\To$ onto the Cuntz algebra $\Oo_n$, and can thus be thought of
as states of $\Oo_n$.  This characterization of essentiality can be extended
to include the case $n=\infty$, even though $\Oo_\infty$ is simple.
The idea is as follows.  When $n=\infty$, Proposition~\ref{prop:isomorphism}
characterizes the state space of a certain concrete \cstar algebra $\Uu$
which contains a copy of $\Oo_\infty$; $\Uu$ also contains
the compact operators $\Kk$.  We give a canonical procedure for extending
states from $\Oo_\infty$ to $\Uu$, and show that a state is essential
if and only if its canonical extension is zero on $\Kk$.
Section~2 concludes with Theorem~\ref{theorem:decomposition},
which gives an alternate approach to the singular-essential decomposition.

A method which has been profitably used to study the state space of $\Oo_n$
has been to focus first on states of its even-word subalgebra,
and then on the problem of extending such states to $\Oo_n$
(\cite{evans, ace, laca1, laca2, bjp, bj}).
When $n$ is finite, this even-word subalgebra is a UHF algebra of
type $n^\infty$, and is thus somewhat less complicated and
better understood than $\Oo_n$.  For example, there is a large
supply of states of this algebra readily at hand in the form of
product states; indeed, these states have played a key r\^ole in the
study of UHF algebras (\cite{glimm,powers}).
In \S3 we study the problem of extending a pure periodic product state
to $\Oo_n$.  To include the case $n=\infty$ in a unified way we reformulate
this problem:  we consider instead the even-word subalgebra $\F$ of
the Toeplitz-Cuntz algebra $\To$,
and focus on pure periodic product states of $\F$ which are essential.
Our main result is Theorem~\ref{theorem:extensions},
which parameterizes the space of all extensions of such a state
by probability measures on the circle.   Pure extensions correspond
to point measures, and we have as an immediate corollary that the space
of pure extensions of such a state is precisely a circle.

The author would like to thank M. Laca for the many helpful discussions
during the preparation of this paper, and W. Arveson for the motivation to
pursue this research.

\section{States of $\To$}\label{section:states}

\subsection*{From $\To$ to $\Kk$ and back}\label{subsection:to and fro}
Suppose $1 \le i_i$, \dots, $i_k \le n$.  We call
$\mu = (i_1,\ldots, i_k)$ a {\em multi-index\/} and define
$\abs\mu := k$ and $v_\mu := v_{i_1}\cdots v_{i_k}$;
of course $v_\emptyset := 1$.
The set of all multi-indices will be denoted $\multi$,
and we define $\multi_k := \setspace{\mu\in\multi: \abs\mu = k}$.
With this notation,
$\To = \clsp\setspace{v_\mu v_\nu^*: \mu,\nu\in\multi}$.

When $n$ is finite, the projection $p := 1 - \sum v_iv_i^*$ generates a closed,
two-sided ideal $\Jj_n$ of $\To$ which is isomorphic to the compact operators
on an infinite-dimensional, separable Hilbert space;
indeed, $\setspace{v_\mu p v_\nu^*: \mu, \nu\in\Ww}$
is a self-adjoint system of matrix units for $\Jj_n$ (\cite{cun}).
Let $q := 1 - p$.  Since $q$ is an identity for $C^*(\set{v_i q})$ 
and $(v_j q)^*(v_iq) = \delta_{ij}q$,
the map $v_i \mapsto v_iq$ extends to a \Star endomorphism $\beta'$ of $\To$.
If we then define $\delta' := \id - \beta'$, one checks easily that
$\delta'(v_\mu v_\nu^*) = v_\mu p v_\nu^*$,
and consequently $\delta'(\To)\subset\Jj_n$.

To include the case $n=\infty$ we implement $\delta'$ spatially utilizing
the Fock representation of $\To$ (\cite{evans}).
Technically speaking, the representation we are about to define
is only unitarily equivalent to the Fock representation;
we prefer this version for purely notational reasons.
Let
\[
\Ee := \clsp\setspace{v_i: 1 \le i \le n},
\]
and more generally, let
\[
\Ee_k := \clsp\setspace{v_\mu:\mu\in\multi_k},\qquad k = 0,1,2,\dots,
\]
so that $\Ee = \Ee_1$.
If $f,g\in\Ee_k$, then $g^*f$ is a scalar multiple of the identity,
and the formula $g^*f = \ip fg 1$ defines an inner product
which makes $\Ee_k$ a Hilbert space.  Notice that the Hilbert space norm
on $\Ee_k$ agrees with the norm $\Ee_k$ inherits as a subspace of $\To$,
and that $\setspace{v_\mu: \mu\in\multi_k}$ is an orthonormal basis for $\Ee_k$.
Let
\[
\Fe := \bigoplus_{k=0}^\infty \Ee_k.
\]
By this we mean nothing more than the abstract direct sum of Hilbert spaces;
in particular, the inclusion maps $\Ee_k \hookrightarrow \To$ do {\em not\/}
factor through the canonical injections $\Ee_k \hookrightarrow\Fe$.
We caution the reader that we will think of $\Ee_k$
in three separate ways: as a subspace of the \cstar algebra
$\To$, as a Hilbert space, and as a subspace of $\Fe$.
This is both a notational advantage and a potential cause of confusion.

For each integer $k \ge 0$, left multiplication by $v_i$ is a linear
isometry from $\Ee_k$ to $\Ee_{k+1}$, and together these maps induce
an isometry $l(v_i)$ on $\Fe$.  Similarly, right multiplication
by $v_i$ induces an isometry $r(v_i)$ on $\Fe$.
Since $l(v_j)^*l(v_i) = \delta_{ij}I$, the map $v_i \mapsto l(v_i)$ extends
to a \Star representation $l$ of $\To$ on $\Fe$; we call this the
{\em Fock representation\/}.  The representation which is
more commonly referred to as the Fock representation is unitarily
equivalent to $l$ via the unitary $\Fe\to\oplus_{k=0}^\infty\Ee^{\otimes k}$
determined by
\[
v_{i_1}\dotsm v_{i_k} \mapsto v_{i_1}\otimes \dots \otimes v_{i_k},
\qquad (i_1,\dots, i_k)\in\Ww.
\]
By \cite{evans}, $l$ is faithful and irreducible.
We will study $\To$ in this representation for the remainder of the paper.

For each pair of vectors $f,g\in\Fe$ we will denote by $\rankone fg$ the
rank-one operator $h\mapsto \ip hg f$ on $\Fe$.
Routine calculations show that when $n<\infty$ we have
$l(v_\mu p v_\nu^*) = \rankone{v_\mu}{v_\nu}$,
so the image of the ideal $\Jj_n$ in the Fock representation is $\Kk$,
the compact operators on $\Fe$.

We implement $\delta'$ spatially as follows.
Define a normal \Star endo\-mor\-phism $\beta$ of $\Bfe$ by
\[
\beta(A) := \sum_{i=1}^n r(v_i)A r(v_i)^*,\qquad A\in\Bfe.
\]
When $n$ is infinite, the above series converges in the strong
operator topology.
One easily checks that $\beta$ implements $\beta'$ spatially
when $n$ is finite; i.e. $\beta(l(v_i)) = l(v_iq)$.
Hence $\delta := \id - \beta$ implements $\delta'$ spatially
when $n<\infty$.  Moreover,
\begin{equation}\label{eq:deltagen}
\delta(l(v_\mu v_\nu^*)) = \rankone{v_\mu}{v_\nu}
\end{equation}
holds whether or not $n$ is finite, so we always have
$\delta\circ l(\To)\subseteq\Kk$.

\subsection*{Operators of Bounded Support}
Let $P_k := I - \beta^{k+1}(I)$, the orthogonal projection of $\Fe$
onto $\bigoplus_{i=0}^k \Ee_i$.
Let $\bs_k$ be the von Neumann algebra of all operators $T\in\Bfe$
satisfying $T = P_k T P_k$, and let
\[
\bs := \bigcup_{k=0}^\infty \bs_k,
\]
the algebra of operators on $\Fe$ which have {\em bounded support\/}.
This algebra is $\beta$\ndash invariant; indeed,
$\beta^i(P_k) = P_{k+i} - P_{i-1}$ for $i\ge 1$, $k\ge 0$.
Consequently
\begin{equation}\label{eq:disjoint}
\beta^i(A)\beta^j(B) = 0 \quad\text{if $A,B\in\bs_k$ and $\abs{i-j}\ge k+1$.}
\end{equation}

Since $\beta^k(I) \to 0$ strongly, 
$\delta$ is injective: if $\delta(A) = 0$, then $A = \beta(A)$,
and hence $A = \beta^k(A) = \beta^k(A)\beta^k(I) = A\beta^k(I) \to 0$.
The following proposition shows that $\bs$ is contained in
the range of $\delta$, and establishes a formula for the inverse of
$\delta$ on $\bs$.

\begin{prop} For each $B\in\bs$, the sum $\sum_{i=0}^\infty \beta^i(B)$
converges $\sigma$\ndash weakly to a bounded operator $\lambda(B)$ on $\Fe$.
The map $\lambda:\bs\to\Bfe$ is injective, \Star linear,
and its restriction to $\bs_k$ is normal, completely positive and
has norm $k+1$.  Moreover, $\delta|_{\lambda(\bs)}$ is the inverse of $\lambda$.
\end{prop}

\begin{proof} We follow \cite[Theorem~2.2]{arv4}.  Fix $k\ge 0$.
If $A,B \in\bs_k$, then by \eqref{eq:disjoint}
\[
\beta^{i(k+1)}(A)\beta^{j(k+1)}(B) = \delta_{ij}\beta^{j(k+1)}(AB).
\]
Consequently $\sum_{j=0}^\infty \beta^{j(k+1)}$
is a normal \Star homomorphism from $\bs_k$ to $\Bfe$.
From this it is easy to see that
$\sum_{i=0}^\infty \beta^i(B)
 = \sum_{l=0}^k \sum_{j=0}^\infty \beta^{j(k+1)+l}(B)$
converges strongly to an operator $\lambda_k(B)$ for each $B\in\bs_k$,
and that $\lambda_k$ is normal, completely positive and \Star linear.
The $\lambda_k$'s are clearly coherent, and thus define the desired map
$\lambda$.  Since $(P_k - P_{k-1})\lambda(P_k) = (k+1)(P_k - P_{k-1})$,
the norm of $\lambda_k$ is at least $k+1$; since $\lambda_k$ is the sum
of $k+1$ \Star homomorphisms, we thus have $\norm{\lambda_k} = k+1$.

By the normality of $\beta$ we have
$\beta(\lambda(B)) = \sum_{i = 1}^\infty \beta^i(B)$,
and thus $\delta(\lambda(B)) = B$ for each $B\in\bs$.
Consequently $\lambda$ is injective with inverse $\delta|_{\lambda(\bs)}$.
\end{proof}

\begin{prop} $\lambda(\bs)$ and $\lambda(\bs\cap\Kk)$ are irreducible
\Star subalgebras of $\Bfe$, and $\lambda(\bs\cap\Kk)$
is a dense subalgebra of $l(\To)$.
\end{prop}

\begin{proof}  Since $\lambda$ is \Star linear,
$\lambda(\bs)$ and $\lambda(\bs\cap\Kk)$
are self-adjoint subspaces of $\Bfe$.
By \eqref{eq:deltagen} we have
$\lambda(\rankone{v_\mu}{v_\nu}) = l(v_\mu v_\nu^*)$,
and since $\lambda$ is bounded when restricted to one of the subalgebras
$\bs_k$, this shows that $\lambda(\bs\cap\Kk)$ is a dense subspace of $l(\To)$.
Suppose $A,B\in\lambda(\bs)$; we will
show that $AB\in\lambda(\bs)$.  To begin with, observe that
if $C\in\bs_c$ and $D\in\bs_d$, then
\begin{equation}\label{eq:bounded}
C\lambda(D) = \sum_{k=0}^c C\beta^k(D) \in \bs_{c + d},
\end{equation}
with a similar equation holding for $\lambda(C)D$.
Since $\delta(A),\delta(B)\in\bs$, this implies that both
$\delta(A)B = \delta(A)\lambda(\delta(B))$ and
$A\delta(B)=\lambda(\delta(A))\delta(B)$ are in $\bs$.  Using the identity
\begin{equation}\label{eq:deltaAB}
\delta(AB) = \delta(A)B + A\delta(B) - \delta(A)\delta(B),
\end{equation}
we see that $\delta(AB)\in\bs$.
Thus $AB\in\lambda(\bs)$, so $\lambda(\bs)$ is a \Star algebra.

If in addition $\delta(A)$ and $\delta(B)$ are compact, then by
\eqref{eq:deltaAB}, $\delta(AB)$ is compact as well.
Thus $AB \in \lambda(\bs\cap\Kk)$, so $\lambda(\bs\cap\Kk)$
is also a \Star algebra.
These algebras are irreducible since $\lambda(\bs\cap\Kk)$
is dense in the irreducible algebra $l(\To)$.
\end{proof}

\begin{remark}\label{remark:ideal}
If $n$ is finite then $\bs\subset\Kk$, so that
the algebras $\lambda(\bs\cap\Kk)$ and $\lambda(\bs)$
coincide.  When $n$ is infinite this is not the case.
For arbitrary $n$, \eqref{eq:bounded} shows that $\overline{\bs}$
is an ideal in $\overline{\lambda(\bs)}$.
Moreover, this ideal contains the compacts
since $\bs$ is irreducible and has nontrivial intersection with $\Kk$.
When $n$ is finite $\overline{\bs} = \Kk$.
\end{remark}

\subsection*{Decreasing Positive Linear Functionals}

\begin{definition}\label{defn:decreasing}
Suppose $\omega$ is a linear functional on $\bs$.  We say that $\omega$
is {\em decreasing\/} if $\omega\circ\delta$ is positive on $\bs$;
that is, if
\[
\omega(B^*B) - \omega(\beta(B^*B)) \ge 0, \qquad B\in\bs.
\]
\end{definition}

\begin{prop}\label{prop:decreasing}
Suppose $\omega$ is a linear functional on $\bs$.
Then $\omega\circ\delta$ is positive on $\lambda(\bs)$
iff $\omega$ is positive and decreasing, in which case
$\omega\circ\delta$ extends uniquely to a positive linear functional
$\Delta\omega$ of norm $\omega(P_0)$ on the \cstar algebra
$\overline{\lambda(\bs)}$.
\end{prop}

\begin{proof} Suppose $\omega\circ\delta$ is positive on $\lambda(\bs)$.
Then $\omega$ is decreasing since $\bs\subseteq\lambda(\bs)$,
and $\omega$ is positive since $\omega = \omega\circ\delta\circ\lambda$
and $\lambda$ is positive.

Conversely, suppose $\omega$ is positive and decreasing.
Then $(B_1,B_2) \mapsto \omega(B_2^*B_1)$  defines a positive semi-definite
sesquilinear form on $\bs$, so there is a Hilbert space
$\Hh_\omega$ and a linear map $\Omega:\bs\to\Hh_\omega$ such that
$\Omega(\bs)$ is dense in $\Hh_\omega$
and $\ip{\Omega(B_1)}{\Omega(B_2)} = \omega(B_2^*B_1)$
for every $B_1,B_2\in\bs$.  Since $\omega$ is decreasing,
$\Omega(B)\mapsto\Omega(\beta(B))$ extends uniquely to a
linear contraction $T$ on $\Hh_\omega$.
Define $\tau:\ZZ\to\Bb(\Hh_\omega)$ by
\begin{equation}\label{eq:tau}
\tau(k) := \begin{cases}
  T^k & \text{if $k\ge0$} \\
  T^{*(-k)} & \text{if $k<0$.}
\end{cases}
\end{equation}
Since $T$ is a contraction, $\tau$ is positive definite.

Suppose now that $B\in\bs$.  By \eqref{eq:deltaAB},
\[
\begin{split}
\delta(\lambda(B)^*\lambda(B))
& = \delta(\lambda(B^*))\lambda(B) +\lambda(B^*)\delta(\lambda(B))
  - \delta(\lambda(B^*))\delta(\lambda(B)) \\
& = B^*\lambda(B) + \lambda(B^*)B - B^*B \\
& = \sum_{i=0}^\infty B^*\beta^i(B) + \sum_{j=1}^\infty \beta^j(B^*)B,
\end{split}
\]
where by \eqref{eq:bounded} each of the sums is only finitely non-zero.
Since
\[
\ip{\tau(k)\Omega(B)}{\Omega(B)} = \begin{cases}
  \omega(B^*\beta^k(B)) & \text{if $k\ge0$} \\
  \omega(\beta^{-k}(B^*)B) & \text{if $k<0$,}
\end{cases}
\]
we thus have
\[
\omega\circ\delta(\lambda(B)^*\lambda(B))
 = \sum_{k = -\infty}^\infty \ip{\tau(k)\Omega(B)}{\Omega(B)} \ge 0,
\]
so $\omega\circ\delta$ is positive on $\lambda(\bs)$.

Since $\lambda(P_0) = \sum_{n=0}^\infty \beta^n(I - \beta(I)) = I$,
the algebra $\lambda(\bs)$ unital.
Hence to show that $\omega\circ\delta$ is bounded on $\lambda(\bs)$
it suffices to establish the relation
\begin{equation}\label{eq:limsup}
\limsup_{k\to\infty} \abs{\omega\circ\delta(A^k)}^{1/k} \le \norm A,
  \qquad A\in\lambda(\bs).
\end{equation}
For this, suppose $A\in\lambda(\bs)$ and $k$ is a positive integer.
Then $\delta(A)\in\bs_a$ for some positive integer $a$, and
repeated applications of \eqref{eq:deltaAB} and \eqref{eq:bounded}
give that $\delta(A^k)\in\bs_{ka}$.
Now $P_{ka}$ is the unit in $\bs_{ka}$, so
\[
\abs{\omega\circ\delta(A^k)} \le \omega(P_{ka})\norm{\delta(A^k)}.
\]
Since $\omega$ is decreasing,
\[
\omega(P_{ka}) = \sum_{i=0}^{ka}\omega(\beta^i(P_0))
\le \sum_{i=0}^{ka}\omega(P_0) = (ka+1)\omega(P_0).
\]
These last two inequalities together with the fact that
$\norm\delta \le 2$ give
\[
\abs{\omega\circ\delta(A^k)} \le 2(ka+1)\omega(P_0)\norm A^k,
\]
from which \eqref{eq:limsup} follows immediately.
Finally,
\[
\norm{\Delta\omega} = \omega\circ\delta(I)
  = \omega\circ\delta(\lambda(P_0)) = \omega(P_0).
\]
\end{proof}

\begin{definition} Denote by $\Pb$ be the cone of
decreasing positive linear functionals on $\bs$,
partial-ordered by the relation
\[
\omega_1 \le \omega_2
\quad\text{iff}\quad \omega_2 - \omega_1 \in\Pb.
\]
\end{definition}  

\begin{prop}\label{prop:isomorphism}
The map $\omega\mapsto \Delta\omega$ is an affine order isomorphism of
$\Pb$ onto the positive part of the dual of $\overline{\lambda(\bs)}$.
\end{prop}

\begin{proof}
It is clear from the construction that $\omega\mapsto\Delta\omega$ is affine.
To see that it is surjective, suppose $\rho$ is a positive linear functional
on $\overline{\lambda(\bs)}$.  Let $\omega = \rho\circ\lambda$.
Then $\omega\circ\delta$ agrees with $\rho$ on $\lambda(\bs)$ and hence
is positive.  By Proposition~\ref{prop:decreasing} this implies that
$\omega$ is positive and decreasing, and clearly $\Delta\omega = \rho$.
Since $\omega(B) = \Delta\omega(\lambda(B))$ for each $B\in\bs$,
the map $\omega\mapsto\Delta\omega$ is injective.
Lastly, observe that $\omega_1 \le \omega_2$ iff $\omega_2 - \omega_1\in\Pb$,
which by Proposition~\ref{prop:decreasing} is equivalent to the condition that
$(\omega_2 - \omega_1)\circ\delta$ be positive on $\lambda(\bs)$.
This in turn is obviously equivalent to the condition that
$\Delta\omega_1 \le \Delta\omega_2$.
\end{proof}

When $n$ is finite recall that
$\overline{\lambda(\bs)} = \overline{\lambda(\bs\cap\Kk)} = l(\To)$,
so we have an affine isomorphism
\begin{equation}\label{eq:Dwl}
\omega \mapsto \Delta\omega\circ l
\end{equation}
of $\Pb$ onto the state space of $\To$.  When $n$ is infinite the algebra
$\overline{\lambda(\bs)}$ properly contains $l(\Oo_\infty)$, and consequently
\eqref{eq:Dwl} gives an affine map of $\Pb$ onto the state space of
$\Oo_\infty$ which is not injective.  The following definition is the key
to identifying a subcone of $\Pb$ on which this map is an isomorphism.

\begin{definition} A linear functional $\omega$ on $\bs$ is
{\em locally normal\/} if its restriction to each of the von Neumann
subalgebras $\bs_k$ is normal.  We denote by $\Wb$ the subcone
of $\Pb$ consisting of all decreasing positive linear functionals on $\bs$
which are locally normal.
\end{definition}

We now prove our main theorem, which improves on
Proposition~\ref{prop:isomorphism} in the sense that it characterizes
the state space not just of $\To$ for $n$ finite, but also of $\Oo_\infty$.

\begin{theorem}\label{theorem:TOn isomorphism}
Suppose $\set{v_1,\dots, v_n}$ are the distinguished generating isometries
of the Toeplitz-Cuntz algebra $\To$; we include the case $n=\infty$ by writing
$\Tt\Oo_\infty$ for the Cuntz algebra $\Oo_\infty$.
Let $\Ee\subseteq\To$ be the closed linear span of $\set{v_1,\dots, v_n}$,
let $\Fe$ be full Fock space over $\Ee$,
let $\bs$ be the algebra of operators on $\Fe$ which have bounded support,
and let $\Wb$ be the partially-ordered cone of decreasing locally normal
positive linear functionals on $\bs$.  For each $\omega\in\Wb$
there is a unique positive linear functional $\rho$ on $\To$ which satisfies
\begin{equation}\label{eq:isomorphism}
\rho(v_\mu v_\nu^*) = \omega(\rankone{v_\mu}{v_\nu}),\qquad \mu,\nu\in\multi.
\end{equation}
Moreover, the map $\omega \mapsto \rho$ is an affine order isomorphism of $\Wb$
onto the positive part of the dual of $\To$.
\end{theorem}

\begin{proof} Suppose $\omega\in\Wb$, and let $\rho := \Delta\omega\circ l$,
where $l$ is the Fock representation of $\To$ and $\Delta\omega$ is as
in Proposition~\ref{prop:decreasing}.  Then $\rho$ satisfies
\eqref{eq:isomorphism}, and $\rho$ is uniquely determined by
\eqref{eq:isomorphism} since elements of the form
$v_\mu v_\nu^*$ have dense linear span in $\To$.
If $n$ is finite then every linear functional on $\bs$ is automatically
locally normal, so that $\Wb$ is all of $\Pb$.
Since $\overline{\lambda(\bs)} = l(\To)$,
the theorem thus reduces to Proposition~\ref{prop:isomorphism}.

It remains only to show that when $n$ is infinite,
\begin{equation}\label{eq:Dwl2}
\omega\mapsto \Delta\omega\circ l
\end{equation}
maps $\Wb$ bijectively onto the positive part of the dual of $\Oo_\infty$.
To begin with, suppose $\omega_1, \omega_2\in\Wb$ are such that
$\Delta\omega_1\circ l = \Delta\omega_2\circ l$.
Then $\Delta\omega_1$ and $\Delta\omega_2$ agree on
$l(\Oo_\infty) \supseteq \lambda(\bs\cap\Kk)$,
so for each $K\in\bs\cap\Kk$ we have
\[
\omega_1(K) = \Delta\omega_1(\lambda(K))
 = \Delta\omega_2(\lambda(K)) = \omega_2(K).
\]
Since $\omega_1$ and $\omega_2$ are locally normal, this implies that
$\omega_1 = \omega_2$.  Thus \eqref{eq:Dwl2} is injective.

To show surjectivity, suppose $\rho$ is a positive linear functional
on $\Oo_\infty$.  Define $\omega_0$ on $\bs\cap\Kk$
by $\omega_0 := \rho\circ l^{-1} \circ \lambda$,
and for each $k$ let $T_k$ be the unique
positive trace-class operator in $\bs_k$ such that
$\omega_0(K) = \tr(T_k K)$ for each $K\in\bs_k\cap\Kk$.
The formula
\begin{equation}\label{eq:Omegak}
\omega(B) := \tr(T_k B),\qquad B\in\bs_k,
\end{equation}
gives the unique extension of $\omega_0$ to a positive linear functional
$\omega$ on $\bs$ which is locally normal.
Once we establish that $\omega$ is decreasing, surjectivity
of \eqref{eq:Dwl2} follows immediately:
for each $K\in\bs\cap\Kk$,
\[
\Delta\omega(\lambda(K)) = \omega\circ\delta(\lambda(K))
  = \omega_0(K) = \rho\circ l^{-1}(\lambda(K)),
\]
so that $\rho = \Delta\omega\circ l$.

We will show that $\omega\circ\delta$ is positive on $\lambda(\bs)$;
by Proposition~\ref{prop:decreasing} this implies that $\omega$ is decreasing.
For each $B\in\bs$ define a function $\phi_B:\ZZ\to\CC$ by
\[
\phi_B(k) := \begin{cases}
\omega(B^*\beta^k(B)) & \text{if $k\ge 0$} \\
\omega(\beta^{-k}(B^*)B) & \text{if $k < 0$.}
\end{cases}
\]
As in the proof of Proposition~\ref{prop:decreasing},
\[
\omega\circ\delta(\lambda(B)^*\lambda(B)) = \sum_{k= -\infty}^\infty \phi_B(k),
\qquad B\in\bs,
\]
so it suffices to show that each $\phi_B$ is positive definite.
We will do this by showing that $\phi_K$ is positive definite
for each $K\in\bs\cap\Kk$ and that any $\phi_B$ can be obtained
as a pointwise limit of such functions.

Let $E_k := \beta^k(I) - \beta^{k+1}(I)$, the orthogonal projection
of $\Fe$ onto $\Ee_k$.  For each $z\in\TT$ let
$U_z := \sum_{k=0}^\infty z^k E_k$.
Each of the unitaries $U_z$ is a multiplier of $\bs\cap\Kk$,
and since $\beta(U_z) = \overline z U_z\beta(I)$ we have
\[
\phi_K(k) = z^k\phi_{U_z K}(k),\qquad K\in\bs\cap\Kk, k\in\ZZ.
\]
Let $\hat\phi_K$ denote the Fourier transform of $\phi_K$.
For each $K\in\bs\cap\Kk$ and $z\in\TT$,
\[
\hat\phi_K(z)
 = \sum_{k = -\infty}^\infty \phi_K (k)z^k
 = \sum_{k = -\infty}^\infty \phi_{U_{\overline z} K}(k)
 = \rho(\lambda(U_{\overline z} K)^*\lambda(U_{\overline z} K))
 \ge 0,
\]
so that $\hat\phi_K$ is positive.  By Herglotz' theorem,
this implies that $\phi_K$ is positive definite.

Lastly, suppose $B\in\bs$, say $B\in\bs_m$.  Let $(K_\alpha)$
be a bounded net in $\bs_m \cap \Kk$ which converges to $B$
in the strong operator topology on $\bs_m$.   Then $K_\alpha^* \to B^*$
in the $\sigma$\ndash weak topology, and hence
$\beta^k(K_\alpha^*) \to \beta^k(B^*)$
$\sigma$\ndash weakly for any $k\ge 0$.
Since this latter net is bounded in norm,
$\beta^k(K_\alpha^*)K_\alpha \to \beta^k(B^*)B$
weakly, hence $\sigma$\ndash weakly.
From this it is apparent that $\phi_{K_\alpha}(k) \to \phi_B(k)$
for each $k\in\ZZ$.  Thus $\phi_B$ is positive definite.
\end{proof}

We conclude this section by giving a reformulation of
Theorem~\ref{theorem:TOn isomorphism} in terms of density matrices.
Suppose $\omega$ is a locally normal linear functional on $\bs$.
Then for each positive integer $k$
there is a unique trace-class operator $T_k$ in $\bs_k$ such that
\[
\omega(B) = \tr(T_k B), \qquad B\in\bs_k.
\]
These density operators are coherent in the sense that
$T_k = P_kT_{k+1}P_k$ for each $k$.
Define $\Omega_{ij}: \Ee_j\to\Ee_i$ by
\[
\Omega_{ij} := E_i T_k E_j,
\]
where $E_i$ is the orthogonal projection of $\Fe$ onto $\Ee_i$,
and $k$ is any integer greater than both $i$ and $j$;
we think of $T_k$ as having operator matrix $(\Omega_{ij})_{i,j=0}^k$.
The {\em density matrix\/} of $\omega$ is the infinite operator matrix
$\Omega:=(\Omega_{ij})$.  If $x\in\Ee_j$ and $y\in\Ee_i$,
we may write $\ip{\Omega x}y$ rather than $\ip{\Omega_{ij}x}y$.

\begin{definition} Suppose $\Omega = (\Omega_{ij})$ is an infinite operator
matrix; i.e., $\Omega_{ij}\in\Bb(\Ee_j,\Ee_i)$ for every pair
$i,j$ of nonnegative integers.  For each $k$ let $T_k$ be the operator
in $\bs_k$ determined by
\[
\ip{T_kx}y = \ip{\Omega x}y
\quad x\in\Ee_j,\, y\in\Ee_i,\, 0 \le i,j \le k.
\]
We say that $\Omega$ is {\em positive\/} if each $T_k$ is positive,
and {\em locally trace-class\/} if each $T_k$ is trace-class.
\end{definition}

It is evident that an infinite operator matrix $\Omega$ is the density matrix
of a locally normal linear functional if and only if it is locally
trace-class, and  that a density matrix $\Omega$ is positive if
and only if its associated linear functional $\omega$ is positive.

We now characterize those density matrices which correspond to linear
functionals which are decreasing.  For this, we first need a lemma.

\begin{lemma}\label{lemma:tc}
Suppose $T$ is a trace-class operator on a separable Hilbert space $\Hh$,
and $U_1$, $U_2$, $U_3$, \dots are isometries
on $\Hh$ with mutually orthogonal ranges.  Then
\begin{equation}\label{eq:tc}
\sum_{k=1}^\infty U_k^*TU_k
\end{equation}
converges in trace-class norm (and hence in operator norm as well).
\end{lemma}

\begin{proof}  First suppose $T\ge 0$.  Let $\set{e_i}$ be an orthonormal basis
for $\Hh$.  If $l>m\ge 1$, then by Fubini's theorem
\begin{align*}
\tr\biggl(\sum_{k=1}^l U_k^*TU_k - \sum_{k=1}^m U_k^*TU_k\biggr)
& = \sum_{i=1}^\infty \sum_{k=m+1}^l \ip{U_k^*TU_ke_i}{e_i} \\
& = \sum_{k=m+1}^l \sum_{i=1}^\infty \ip{TU_ke_i}{U_ke_i}.
\end{align*}
But
\[
\sum_{k=1}^\infty \sum_{i=1}^\infty \ip{TU_ke_i}{U_ke_i} \le \tr T <\infty
\]
so the sequence $(\sum_{k=1}^m U_k^*TU_k)_{m=1}^\infty$
is Cauchy in the trace-class norm.
Since the algebra of trace-class operators is complete in this norm,
the infinite sum \eqref{eq:tc} converges as claimed.
Since every trace-class operator can be written
as a linear combination of four positive trace-class operators,
\eqref{eq:tc} converges in trace-class norm for every trace-class
operator $T$.
\end{proof}

\begin{definition} Suppose $T$ is a trace-class operator on $\Fe$.
The {\em slice\/} of $T$ is the operator
\[
\slice T := \sum_{k=1}^n r(v_k)^*Tr(v_k),
\]
where $r(v_k)$ is right creation by $v_k$ on $\Fe$.  Note that
when $n=\infty$ the sum converges to a trace-class operator
(Lemma~\ref{lemma:tc}).

If $\Omega = (\Omega_{ij})$ is a locally trace-class
operator matrix, we denote by $\slice\Omega$ the locally trace-class
operator matrix $(\slice\Omega_{i+1,j+1})$.
\end{definition}

\begin{remark} Suppose $S\in \Bb(\Ee_j,\Ee_i)$ and $A\in\Bb(\Ee)$
are trace-class operators such that $T(xy) = (Sx)(Ay)$ for
$x\in\Ee_j$ and $y\in\Ee$.  (The unitary $x\otimes y \mapsto xy$
transforms $T$ into the tensor product $S\otimes A$.)  Then
\begin{multline*}
(\slice T)x = \sum_{k=1}^n r(v_k)^* T r(v_k)x
= \sum_{k=1}^n r(v_k)^* T(xv_k) \\
= \sum_{k=1}^n r(v_k)^* (Sx)(Av_k)
= \sum_{k=1}^n \ip{Av_k}{v_k}Sx
= \tr(A)Sx,
\end{multline*}
so slicing has the effect of taking the trace in the last variable.
\end{remark}

\begin{lemma}\label{lemma:slice}
Suppose $\omega$ is a locally-normal linear functional on $\bs$ with
density matrix $\Omega$.  Then $\omega\circ\beta$ has density matrix
$\slice\Omega$.
\end{lemma}

\begin{proof} Let $T_k$ be the trace-class operator in $\bs_k$ such that
$\omega(B) = \tr(T_kB)$ for $B\in\bs_k$, so that $T_k$
has operator matrix $(\Omega_{ij})_{i,j=0}^k$.
Suppose $B\in\bs_k$.  Then $\beta(B)\in\bs_{k+1}$, and thus
\begin{multline*}
\omega\circ\beta(B)
= \tr(T_{k+1}\beta(B))
= \sum_{i=1}^\infty \tr(T_{k+1} r(v_i)Br(v_i)^*) \\
= \sum_{i=1}^\infty \tr(r(v_i)^*T_{k+1} r(v_i)B)
= \tr((\slice T_{k+1})B).
\end{multline*}
Thus $\omega\circ\beta$ has density matrix $\slice\Omega$.
\end{proof}

We now give our reformulation of Theorem~\ref{theorem:TOn isomorphism}.

\begin{theorem}\label{theorem:TOn isomorphism2}
Suppose $\set{v_1,\dots, v_n}$ are the distinguished generating iso\-metries
of the Toeplitz-Cuntz algebra $\To$; we include the case $n=\infty$ by writing
$\Tt\Oo_\infty$ for the Cuntz algebra $\Oo_\infty$.
Let $\Ee\subseteq\To$ be the closed linear span of $\set{v_1,\dots, v_n}$,
and let $\Fe$ be full Fock space over $\Ee$.
Suppose $\Omega$ is a positive locally trace-class operator matrix
on $\Fe$ which satisfies $\slice\Omega\le\Omega$.
Then there is a unique positive linear functional $\rho$ on $\To$ which
satisfies
\[
\rho(v_\mu v_\nu^*) = \ip{\Omega v_\mu}{v_\nu},\qquad \mu,\nu\in\multi.
\]
Moreover, the map $\Omega \mapsto \rho$ is an affine order isomorphism of such
operator matrices onto the positive part of the dual of $\To$.
\end{theorem}

\begin{proof} The equation
$\omega(\rankone{v_\mu}{v_\nu}) = \ip{\Omega v_\mu}{v_\nu}$
establishes an affine order isomorphism $\omega\mapsto\Omega$ between
$\Wb$ and positive locally trace-class operator matrices
on $\Fe$ which satisfies $\slice\Omega\le\Omega$, so the theorem
follows immediately from Theorem~\ref{theorem:TOn isomorphism}.
\end{proof}

\section{Singular and Essential States}\label{section:sing-ess}

A state $\rho$ of $\To$ with GNS representation $\pi:\To\to\bh$
is said to be {\em essential\/} if $\sum \pi(v_iv_i^*)$ 
is the identity operator on $\Hh$, and {\em singular\/}
if $\sum_{\mu\in\multi_k} \pi(v_\mu v_\mu^*)$ decreases strongly to zero in $k$.
When $n$ is finite $\To$ has a unique ideal $\Jj_n$,
and it is not hard to show that essential states of $\To$
are precisely those which are singular with respect to this ideal;
similarly, singular states are $\Jj_n$\ndash essential.
As a result, every state of $\To$ has a unique decomposition
into essential and singular components, a result which was
generalized to the case $n=\infty$ in \cite{laca1}.

We can view the singular/essential decomposition of a state $\rho$
of $\Oo_\infty$ as the decomposition with respect to an ideal as follows.
Let $\omega$ be the unique decreasing locally normal positive linear functional
on $\bs$ such that $\rho = \Delta\omega\circ l$,
as in Theorem~\ref{theorem:TOn isomorphism}.  The \cstar algebra
$\overline{\lambda(\bs)}$ contains the ideal $\Kk$
of compact operators on $\Fe$ (Remark~\ref{remark:ideal}),
so we can decompose the functional
$\Delta\omega$ of Proposition~\ref{prop:isomorphism} with respect to this
ideal.  Restricting to $\overline{\lambda(\bs\cap\Kk)} = l(\Oo_\infty)$
and pulling back to $\Oo_\infty$ gives the singular/essential decomposition
of $\rho$; this will follow from Proposition~\ref{prop:sing-ess}~(1d) and (2d).

Again allowing $n$ to be either finite or infinite, we follow
\cite{laca1} and define a positive linear functional $\alpha^*\rho$ by
\[
\alpha^*\rho(x) := \sum_{i=1}^n \rho(v_i x v_i^*),\qquad x\in\To.
\]
In \cite[Corollary~2.9]{laca1}, Laca characterized singular and essential
states of $\To$ using the monotonically nonincreasing sequence
$(\norm{\alpha^{*k}\rho})_{k=1}^\infty$:  $\rho$ is essential iff this sequence is constant
and singular iff it converges to zero.

Let $\omega$ be such that $\rho = \Delta\omega\circ l$, as in
Theorem~\ref{theorem:TOn isomorphism}.  Then $\omega \circ \beta$
is a locally normal positive linear functional on $\bs$ which is decreasing
since $(\omega\circ\beta)\circ\delta = (\omega\circ\delta)\circ\beta$
is positive on $\bs$.  Define
\[
\beta^*\rho := \Delta(\omega\circ\beta)\circ l.
\]
If $\mu, \nu \in\multi$, then
\begin{multline}
\beta^*\rho(v_\mu v_\nu^*)
 = \omega\circ\beta(\rankone{v_\mu}{v_\nu}) 
 = \omega\biggl(\sum_{i=1}^n r(v_i)(\rankone{v_\mu}{v_\nu})r(v_i)^*\biggr) \\
 = \sum_{i=1}^n \omega(\rankone{v_\mu v_i}{v_\nu v_i})
 = \sum_{i=1}^n \rho(v_\mu v_i v_i^* v_\nu^*).
\end{multline}
In particular
$\norm{\beta^{*k}\rho}
= \beta^{*k}\rho(1)
= \sum_{\mu\in\multi_k} \rho(v_\mu v_\mu^*)
= \norm{\alpha^{*k}\rho}$,
so Laca's characterization can be stated in terms of $\beta^*$.
What is not immediately apparent is that essentiality is equivalent to
$\beta^*$\ndash invariance.

We remind the reader of the notation $E_k := \beta^k(I) - \beta^{k+1}(I)$,
the orthogonal projection of $\Fe$ onto $\Ee_k$.

\begin{prop}\label{prop:sing-ess}
Suppose $\rho$ is a positive linear functional on $\To$
and $\omega$ is the unique decreasing locally normal positive
linear functional on $\bs$ such that $\rho = \Delta\omega \circ l$.
Let $\Omega$ be the density matrix of $\omega$.
Statements \textup{(1a)--(1f)} below are equivalent, as are statements
\textup{(2a)--(2d)}:
\begin{quote}
\hfil
\begin{tabular}{l}
 \textup{(1a)} $\rho$ is essential \\
 \textup{(1b)} $\omega(E_k)$ is constant in $k$ \\
 \textup{(1c)} $\omega\circ\delta(B) = 0$ for each $B\in\bs$ \\
 \textup{(1d)} $\Delta\omega$ is $\Kk$\ndash singular \\
 \textup{(1e)} $\slice\Omega = \Omega$ \\
 \textup{(1f)} $\rho = \beta^*\rho$
\end{tabular}
\hfil
\begin{tabular}{l}
 \textup{(2a)} $\rho$ is singular \\
 \textup{(2b)} $\lim \omega(E_k) = 0$ \\
 \textup{(2c)} $\rho$ is normal in the \\
 \phantom{\textup{(2c)}}Fock representation \\
 \textup{(2d)} $\Delta\omega$ is $\Kk$\ndash essential.
\end{tabular}
\hfil
\end{quote}
\end{prop}

\begin{proof} Since $\setspace{v_\mu:\mu\in\multi_k}$ is an orthonormal
basis for $\Ee_k$ and $\omega$ is locally normal,
\[
\omega(E_k)
= \sum_{\mu\in\multi_k} \omega(\rankone{v_\mu}{v_\mu})
= \sum_{\mu\in\multi_k} \rho(v_\mu v_\mu^*)
= \norm{\alpha^{*k}\rho}.
\]
Thus (1a) $\iff$ (1b) and (2a) $\iff$ (2b) follow from
\cite[Corollary~2.9]{laca1}.

(1b) $\iff$ (1c):  Since $\delta(E_k) = E_k - E_{k+1}$,
(1c) $\Ra$ (1b) is immediate.
For the converse, simply observe that $\omega\circ\delta|_{\bs}$
is a positive linear
functional whose restriction to the \cstar algebra $\bs_k$ has norm
$\omega\circ\delta(P_k) = \omega\circ\delta(\sum_{i=0}^k E_i) = 0$.

(1c) $\iff$ (1d):  Since $\bs \subseteq \lambda(\bs)$
and $\Delta\omega = \omega\circ\delta$
on $\lambda(\bs)$, (1c) implies that $\Delta\omega$ vanishes on $\bs$,
and hence on $\overline{\bs}$.  Since $\Kk\subseteq\overline{\bs}$,
this gives (1d).
Conversely, if $\Delta\omega(K) = 0$ for each $K\in\Kk$,
then $\omega\circ\delta(K) = 0$  for each $K\in\bs\cap\Kk$.
Since $\delta$ is $\sigma$\ndash weakly continuous and maps
$\bs_k$ into $\bs_{k+1}$, (1c) follows from local normality of $\omega$.

(1c) $\iff$ (1e):  By Lemma~\ref{lemma:slice}, 
$\omega = \omega\circ\beta$ iff $\slice\Omega = \Omega$.

(1c) $\iff$ (1f):  By Proposition~\ref{prop:isomorphism}, 
$\omega = \omega\circ\beta$ iff $\rho = \beta^*\rho$.

(2a) $\Ra$ (2c):  This follows from \cite[Theorem~2.11]{laca1}.

(2c) $\Ra$ (2b):  Suppose $\rho = \phi \circ l$ for some $\phi\in\Bfe_*$.
Then $\phi$ and $\Delta\omega$ agree on $l(\To)$, so $\phi \circ \lambda$
and $\omega$ agree on $\bs\cap\Kk$.
Fix $k$, and let $\set{P_\alpha}$
be a net of finite rank projections which increases to $E_k$.
By the normality of $\phi$ and local normality of $\lambda$,
\[
\omega(E_k)
 = \lim \omega(P_\alpha)
 = \lim \phi\circ\lambda(P_\alpha)
 = \phi\circ\lambda(E_k)
 = \phi(\beta^k(I)),
\]
which decreases to zero.

(2c) $\iff$ (2d):  Suppose again that $\rho = \phi \circ l$
for some $\phi\in\Bfe_*$.  Fix $B\in\bs$, say $B\in\bs_m$,
and let $\set{K_\alpha}$ be a net of compact operators
in $\bs_m$ which converge $\sigma$\ndash weakly to $B$.  Then
\[
\Delta\omega(\lambda(B))
 = \omega(B)
 = \lim\omega(K_\alpha)
 = \lim\phi\circ\lambda(K_\alpha)
 = \phi(\lambda(B)),
\]
so that $\phi$ extends $\Delta\omega$ as well.  The converse is immediate.
\end{proof}

\begin{theorem}\label{theorem:decomposition}
Suppose $\rho$ is a positive linear functional on $\To$,
and let $\omega$ be the unique decreasing locally normal positive linear
functional on $\bs$ such that $\rho = \Delta\omega \circ l$.
Then $\omega\circ\delta|_{\bs}$ extends uniquely to a normal positive linear
functional $\phi$ on $\Bfe$, and the singular part of $\rho$ is $\phi\circ l$.
\end{theorem}

\begin{proof} Suppose $\omega$ is a decreasing locally normal positive linear
functional on $\bs$.  Then $\omega\circ\delta|_{\bs}$
is a locally normal positive linear functional,
and $\omega\circ\delta(P_k) = \omega(E_0) - \omega(E_{k+1}) \le \omega(E_0)$
for every $k\ge0$.  It follows that $\omega\circ\delta$ is bounded on $\bs$:
if $B\in\bs_k$, then 
$\abs{\omega\circ\delta(B)} \le \omega\circ\delta(P_k)\norm B
 \le \omega(E_0)\norm B$.
Thus $\omega\circ\delta$ extends uniquely to a positive linear functional 
$\psi$ on $\overline{\bs}$.  By Remark~\ref{remark:ideal}
$\Kk\subset\overline{\bs}$,
so there is a unique $\phi\in\Bfe_*$ which coincides with $\psi$ on $\Kk$.
But then $\phi(K) = \omega\circ\delta(K)$ for every $K\in\bs\cap\Kk$,
which by local normality implies that $\phi$ extends $\omega\circ\delta|_{\bs}$.

Let $\rho_s := \phi\circ l$.  By Proposition~\ref{prop:sing-ess}~(2c),
$\rho_s$ is singular.  Let $\omega_s$ be the unique decreasing locally
normal positive linear functional on $\bs$ such that
$\rho_s = \Delta\omega_s \circ l$, and let $\omega_e$ be the locally normal
linear functional $\omega - \omega_s$.  By Theorem~\ref{theorem:TOn isomorphism}
and Proposition~\ref{prop:sing-ess}~(1c), the proof will be complete once
we establish that $\omega_e$ is positive and $\omega_e\circ\delta|_{\bs} = 0$.

For every $K\in\bs\cap\Kk$ we have
$\omega_e(K) = \omega(K) - \omega_s(K) = \omega(K) - \phi(\lambda(K))$,
so by local normality we have $\omega_e = \omega - \phi\circ\lambda$. 
Consequently $\omega_e\circ\delta|_{\bs} = 0$.

To show that $\omega_e$ is positive, we fix a positive integer $m$
and show that the bounded linear functional $\omega_e|_{\bs_m}$
achieves its norm at $P_m$, the identity element of the \cstar algebra $\bs_m$.

For every $k\ge0$ we have $\omega_e = \omega_e\circ\beta^k$, so
\begin{align*}
\norm{\omega_e|_{\bs_m}}
& = \norm{\omega_e\circ\beta^k|_{\bs_m}} \\
& \le \norm{\omega\circ\beta^k|_{\bs_m}}
  + \norm{\omega_s\circ\beta^k|_{\bs_m}} \\
& = \omega\circ\beta^k(P_m) + \omega_s\circ\beta^k(P_m) \\
& = \sum_{i=0}^m (\omega(E_{k+i}) + \omega_s(E_{k+i})) \\
& \le (m+1)(\omega(E_k) + \omega_s(E_k)),
\end{align*} 
since $\omega(E_k)$ and $\omega_s(E_k)$ are monotonically nonincreasing in $k$.
By Proposition~\ref{prop:sing-ess}~(2b), $\lim_{k\to\infty} \omega_s(E_k) = 0$,
so $\norm{\omega_e|_{\bs_m}} \le (m+1)\lim_{k\to\infty}\omega(E_k)$.
On the other hand,
\begin{align*}
\omega_e(P_m)
& = \omega_e\circ\beta^k(P_m) \\
& = \sum_{i=0}^m (\omega(E_{k+i}) - \omega_s(E_{k+i})) \\
& \ge (m+1) (\omega(E_{k+m}) - \omega_s(E_k)),
\end{align*}
so $\omega_e(P_m) \ge (m+1)\lim_{k\to\infty} \omega(E_k)$.
Thus $\omega_e$ is positive.
\end{proof}

\section{Extending Product States}\label{section:extensions}

For each $\lambda \in \TT$ the isometries
$\setspace{\lambda v_i: 1 \le i \le n}$
satisfy the relations $(\lambda v_j)^*(\lambda v_i) = \delta_{ij} 1$,
and hence there is a \Star endomorphism $\gamma_\lambda$ of $\To$
such that $\gamma_\lambda(v_i) = \lambda v_i$.  Each $\gamma_\lambda$
is actually an automorphism since $\gamma_\lambda \circ \gamma_{\lambda^{-1}}$
is the identity on $\To$; in fact $\gamma$ is a continuous automorphic action
of $\TT$ on $\To$, called the {\em gauge\/} action.  Denote by $\F$ the
fixed-point algebra of this action, and let $\Phi$ denote the
canonical conditional expectation of $\To$ onto $\F$; that is,
\[
\Phi(x) := \int_{\TT} \gamma_\lambda(x)\,dm(\lambda), \qquad x\in \To,
\]
where $m$ is normalized Haar measure.  In terms of generating monomials,
\[
\F = \clsp\setspace{v_\mu v_\nu^*: \abs\mu = \abs\nu}
\quad\text{and}\quad
\Phi(v_\mu v_\nu^*)
  = \begin{cases}
  v_\mu v_\nu^* & \text{if $\abs\mu = \abs\nu$} \\
  0 & \text{otherwise.} \end{cases}
\]

We are interested in {\em product states\/} of $\F$.
To explain what we mean by this, let
$\tilde\Kk := \Kk(\Ee) \times \CC$,
endowed with the structure of a unital \Star algebra via
$(A, \lambda)(B, \mu) = (AB + \lambda B + \mu A, \lambda\mu)$
and $(A, \lambda)^* = (A^*, \overline \lambda)$.
When $\Ee$ is infinite-dimensional $\tilde\Kk$ is \Star isomorphic
to the concrete \cstar algebra $\Kk(\Ee) + \CC I$,
but when $\dim\Ee<\infty$ this is not the case.
Nevertheless, it is not difficult to show that there is a unique
\cstar norm on $\tilde\Kk$.
The map
\begin{multline*}
v_{i_1}\dotsm v_{i_m}v_{j_m}^*\dotsm v_{j_1}^*
\mapsto \\
(\rankone{v_{i_1}}{v_{j_1}}, 0) \otimes \dots \otimes
(\rankone{v_{i_m}}{v_{j_m}}, 0)
  \otimes (0, 1) \otimes (0, 1) \otimes \dotsm
\end{multline*}
embeds $\F$ in the infinite tensor product $\tilde\Kk^{\otimes \infty}$.
If $(\rho_k)_{k=1}^\infty$ is a sequence of states of $\Kk(\Ee)$,
so that $\tilde\rho_k(K, \lambda) = \rho_k(K) + \lambda$
defines a sequence of states of $\tilde\Kk$, we call the restriction
of the product state $\otimes_{k=1}^\infty \tilde\rho_k$ to $\F$
a {\em product state\/} of $\F$.

Now suppose $(e_k)_{k=1}^\infty$ is a sequence of unit vectors in $\Ee$.
For each $k$, let $\rho_k$ denote the vector state of $\Kk(\Ee)$
corresponding to $e_k$,
and let $\rho$ denote the corresponding product state of $\F$.
It is evident that $\rho$ is pure and determined by
\begin{equation}\label{eq:rho}
\rho(v_{i_1}\dotsm v_{i_m}v_{j_m}^*\dotsm v_{j_1}^*)
= \ip{v_{i_1}}{e_1}\dotsm\ip{v_{i_m}}{e_m}
  \ip{e_m}{v_{j_m}}\dotsm\ip{e_1}{v_{j_1}}.
\end{equation}
The remainder of this paper is devoted to classifying all extensions
to $\To$ of such a state.

One can always extend $\rho$ by precomposing with $\Phi$; the resulting
extension $\rho \circ \Phi$ is called the {\em gauge-invariant\/} extension.
The most extreme situation is when this extension is pure, in which case
it is the unique state extending $\rho$.
By \cite[Theorem~4.3]{laca2}, this occurs precisely when the sequence $(e_k)$
is {\em aperiodic\/} in the sense that the series
\begin{equation}\label{eq:period}
\sum_{i=1}^\infty (1 - \abs{\ip{e_i}{e_{i+p}}})
\end{equation}
diverges for each positive integer $p$.
In all other cases we say that $(e_k)$ is {\em periodic\/},
and call the smallest positive integer $p$
for which the series in \eqref{eq:period} converges the {\em period\/}
of $(e_k)$.

Suppose then that $(e_k)$ has finite period $p$.
Notice from \eqref{eq:rho} that if we multiply each of the
vectors $e_k$ by a complex number of modulus one, we obtain
a sequence which gives rise to the same product state $\rho$.
Consequently, we are free to rephase so that $\ip{e_i}{e_{i+p}}$
is always real and nonnegative.

\begin{theorem}\label{theorem:extensions}
Suppose $\set{v_1,\dots, v_n}$ are the distinguished generating isometries
of the Toeplitz-Cuntz algebra $\To$; we include the case $n=\infty$ by writing
$\Tt\Oo_\infty$ for the Cuntz algebra $\Oo_\infty$.
Let $\Ee\subseteq\To$ be the closed linear span of $\set{v_1,\dots, v_n}$,
and suppose $(e_k)$ is a sequence of unit vectors in $\Ee$ which is periodic
with finite period $p \ge 1$, and for which $\ip{e_i}{e_{i+p}}$ is always
nonnegative.  Let $\rho$ be the corresponding pure product state of $\F$
determined by \eqref{eq:rho}.
There is an affine isomorphism $\sigma \mapsto \rho_\sigma$
from $P(\TT)$, the space of Borel probability measures on the circle $\TT$,
to the space of all states of $\To$ which extend $\rho$,
given by
\begin{equation}\label{eq:rhosigma}
\rho_\sigma(v_{i_1}\dotsm v_{i_k}v_{j_l}^*\dotsm v_{j_1}^*)
= \lambda_{k,l}
  \ip{v_{i_1}}{e_1} \dotsm \ip{v_{i_k}}{e_k}
  \ip{e_l}{v_{j_l}} \dotsm \ip{e_1}{v_{j_1}},
\end{equation}
where
\begin{equation}\label{eq:coeff}
\lambda_{k,l}
  = \begin{cases}
      \hat\sigma\left(\frac{k-l}p \right)
      \prod_{i=1}^\infty \ip{e_{l+i}}{e_{k+i}}
        & \text{if $k-l \in p\ZZ$} \\
      0 & \text{otherwise,}
\end{cases}
\end{equation}
and $\hat\sigma$ is the Fourier transform of $\sigma$.

\end{theorem}

\begin{proof}
For each $k \ge 0$ let $\be_k := e_1 \dotsm e_k \in \To$,
where of course $\be_0 := 1$.
Note that $\be_k \in \Ee_k = \clsp\setspace{v_\mu: \mu\in\multi_k}$.
Equation \eqref{eq:rho} can now be written more tersely as
\[
\rho(v_\mu v_\nu^*) = \ip{v_\mu}{\be_{\abs\mu}}\ip{\be_{\abs\nu}}{v_\nu},
\qquad \mu,\nu\in\multi, \, \abs\mu = \abs\nu,
\]
which in turn extends by linearity and continuity to
\[
\rho(xy^*) = \ip x{\be_k} \ip {\be_k}y,\qquad x,y\in\Ee_k.
\]
Similarly, \eqref{eq:rhosigma} can be written as
\begin{equation}\label{eq:rhosigma2}
\rho_\sigma(xy^*)
 = \lambda_{k,l} \ip x{\be_k}\ip{\be_l}y,
\qquad x\in\Ee_k, \, y\in\Ee_l.
\end{equation}

Suppose now that $\tilde\rho$ is a state of $\To$ which extends $\rho$.
By the Schwarz inequality,
\begin{equation}\label{eq:cs}
\abs{\tilde\rho(xy^*)}
 \le \rho(xx^*)^{1/2} \rho(yy^*)^{1/2}
 = \abs{\ip x{\be_k}\ip y{\be_l}},
\qquad x\in\Ee_k, \, y\in\Ee_l.
\end{equation}
Let $x_1 := \ip x{\be_k} \be_k$,    $x_2 := x - x_1$,
    $y_1 := \ip y{\be_l} \be_l$ and $y_2 := y - y_1$.
Then
\[
\tilde\rho(xy^*)
= \tilde\rho(x_1 y_1^*) + \tilde\rho(x_1 y_2^*)
  + \tilde\rho(x_2 y_1^*) + \tilde\rho(x_2 y_2^*)
= \tilde\rho(x_1 y_1^*);
\]
that is,
\begin{equation}\label{eq:lambda}
\tilde\rho(xy^*)
= \tilde\rho(\be_k\be_l^*) \ip x{\be_k} \ip{\be_l}y,
\qquad x\in\Ee_k, y\in\Ee_l.
\end{equation}
For each $k,l \ge 0$ define
\[
\lambda_{k,l} := \tilde\rho(\be_k\be_l^*).
\]
Comparing \eqref{eq:rhosigma2} and \eqref{eq:lambda},
it is apparent that we must exhibit a Borel probability measure $\sigma$
such that \eqref{eq:coeff} is satisfied.

For each positive integer $k$
\[
\sum_{\mu\in\multi_k} \tilde\rho(v_\mu v_\mu^*)
 = \sum_{\mu\in\multi_k}  \abs{\ip{v_\mu}{\be_k}}^2
 = \norm{\be_k}^2 = 1,
\]
so by \cite[Corollary~2.9]{laca1} $\tilde\rho$ is essential.
By Proposition~\ref{prop:sing-ess}~(1f), this implies that
$\tilde\rho = \beta^*\tilde\rho$.  In particular, for each $k,l\ge 0$
\begin{align*}
\lambda_{k,l}
& = \tilde\rho(\be_k \be_l^*)
  = \beta^*\tilde\rho(\be_k \be_l^*)
  = \sum_{i=1}^n \tilde\rho(\be_k v_i v_i^* \be_l^*) \\
& = \sum_{i=1}^n \tilde\rho(\be_{k+1} \be_{l+1}^*)
                 \ip{\be_k v_i}{\be_{k+1}}\ip{\be_{l+1}}{\be_l v_i}
    \qquad\text{(by \eqref{eq:lambda})} \\
& = \sum_{i=1}^n \lambda_{k+1,l+1}
                 \ip{v_i}{e_{k+1}}\ip{e_{l+1}}{v_i} \\
& = \lambda_{k+1,l+1} \ip{e_{l+1}}{e_{k+1}},
\end{align*}
and by induction
\begin{equation}\label{eq:recursion}
\lambda_{k,l}
= \lambda_{k+j,l+j} \prod_{i=1}^j \ip{e_{l+i}}{e_{k+i}}
\qquad \forall j.
\end{equation}

Suppose now that $p$ divides $\abs{k - l}$.
Since we have phased the sequence $(e_i)$ so that
$\ip{e_i}{e_{i+p}}$ is always nonnegative, the assumption
that $(e_i)$ has period $p$ means that
$\sum \abs{1 - \ip{e_i}{e_{i+p}}} < \infty$.
By \cite[Proposition~1.2]{gui}, it follows that
$\sum \abs{1 - \ip{e_i}{e_{i+m}}} < \infty$
whenever $p$ divides $m$.  In particular we have
\[
\sum_{i=1}^\infty \abs{1 - \ip{e_{l+i}}{e_{k+i}}} < \infty,
\]
which implies that there is a positive integer $i_0$ such that
\[
\lim_{j\to\infty} \prod_{i = i_0}^{i_0 + j} \ip{e_{l+i}}{e_{k+i}}
\]
exists and is nonzero.
Together with \eqref{eq:recursion}, this shows that
$\lim_{m \to \infty} \lambda_{k+m,l+m}$ exists; indeed
\[
\lim_{m \to \infty} \lambda_{k+m,l+m}
 = \lambda_{k + i_0 - 1, l + i_0 - 1}
   \biggl( \prod_{i = i_0}^\infty \ip{e_{l+i}}{e_{k+i}} \biggr)^{-1}.
\]
Since this limit depends only on the quantity $k-l$,
we can define a function $\tau:\ZZ \to \CC$
by $\tau_{a-b} := \lim_{m \to \infty} \lambda_{ap+m,bp+m}$ for $a,b \ge 0$.

We claim that $\tau$ is positive definite; that is, we claim that for any
finite collection $z_0$, \dots, $z_m$ of complex numbers the sum
$\sum_{a,b = 0}^m z_a\overline{z_b} \tau_{a-b}$ is real and nonnegative.
To see this, define a sequence $(w_i)$ in $\To$ by
$w_i := \sum_{a = 0}^m z_a \be_{(a+i)p}$.
Then
\begin{multline*}
\sum_{a,b = 0}^m z_a \overline{z_b} \tau_{a-b}
 = \lim_{i\to\infty} \sum_{a,b = 0}^m
                     \lambda_{(a+i)p,(b+i)p} z_a\overline{z_b} \\
 = \lim_{i\to\infty} \sum_{a,b = 0}^m
                     \tilde\rho(\be_{(a+i)p}\be_{(b+i)p}^*) z_a\overline{z_b}
 = \lim_{i\to\infty} \tilde\rho(w_i w_i^*)
\ge 0,
\end{multline*}
as claimed.

By Herglotz' theorem there is Borel probability measure $\sigma$ on $\TT$
such that $\tau = \hat\sigma$.
We claim that \eqref{eq:coeff} is satisfied.  The case $k-l\in p\ZZ$
follows immediately from \eqref{eq:recursion} by letting $j\to\infty$.
If $p$ does not divide $k-l$, then by \cite[Proposition~1.2]{gui}
the series $\sum_{i=1}^\infty (1 - \abs{\ip{e_{l+i}}{e_{k+i}}})$
diverges, so that the infinite product
$\prod_{i=1}^\infty \abs{\ip{e_{l+i}}{e_{k+i}}}$ diverges as well;
in particular,
\[
\lim_{j\to\infty} \prod_{i=1}^j \ip{e_{l+i}}{e_{k+i}} = 0.
\]
Since by \eqref{eq:cs} we have $\abs{\lambda_{k+j,l+j}} \le 1$ for each $j$,
it follows from \eqref{eq:recursion} that $\lambda_{k,l} = 0$.
This completes the proof that every extension of $\rho$
is of the form $\rho_\sigma$.

Conversely, suppose $\sigma$ is a Borel probability measure on $\TT$.
Define coefficients $\lambda_{k,l}$ as in \eqref{eq:coeff},
and, resuming the notation and terminology of \S\ref{section:states},
define a locally normal linear functional $\omega$ on $\bs$ by
\begin{equation}\label{eq:omegaB}
\omega(B) := \sum_{k,l = 0}^m \lambda_{k,l} \ip{B\be_l}{\be_k},\qquad B\in\bs_m.
\end{equation}
We claim that $\omega$ is positive and decreasing, and that the functional
$\Delta\omega\circ l$ of Proposition~\ref{prop:decreasing}
is the desired state $\rho_\sigma$ satisfying \eqref{eq:rhosigma2}.

For $c = 0$, $1$, \dots, $p-1$, let $\Hh_c$ be the Hilbert space
inductive limit of the isometric inclusions $\Ee_m \hookrightarrow \Ee_{m+1}$
determined by
\[
x_1\dotsm x_m \mapsto x_1 \dotsm x_m e_{m+c+1},\qquad x_i\in\Ee.
\]
Modulo the isomorphisms
$x_1\dotsm x_m\in\Ee_m \mapsto
 x_1 \otimes \dots \otimes x_m \in \Ee^{\otimes m}$,
$\Hh_c$ is just the infinite tensor product $\Ee^{\otimes \infty}$ with
canonical unit vector $e_{c+1} \otimes e_{c+2} \otimes e_{c+3} \otimes \dotsm$
introduced in \cite{vonNeumann}.
Consequently \cite[Proposition~1.1]{gui} applies:
if $(f_i)$ is a sequence of unit vectors
in $\Ee$ such that $\sum_{i=1}^\infty \abs{1 - \ip{e_{c+i}}{f_i}} < \infty$,
then $f_1 f_2 f_3 \dotsm$ is a unit vector in $\Hh_c$.
In particular, for each $a \ge 0$ we can define a vector
$f_{c,a} \in \Hh_c$ by
\[
f_{c,a} := e_{ap + c + 1} e_{ap + c + 2} e_{ap + c + 3} \dotsm.
\]
By \eqref{eq:coeff},
$\lambda_{ap + c, bp + c} = \hat\sigma(a - b)\ip{f_{c,b}}{f_{c,a}}$.

Suppose now that $B$ is an operator of bounded support on $\Fe$.
Choose $M$ so that $B\in\bs_{Mp + p - 1}$.  Then
\begin{align*}
\omega(B^*B)
& = \sum_{k,l = 0}^{Mp + p - 1} \lambda_{k,l} \ip{B\be_l}{B\be_k} \\
& = \sum_{c=0}^{p-1} \sum_{a,b = 0}^M
    \lambda_{ap + c, bp + c} \ip{B\be_{bp + c}}{B\be_{ap + c}} \\
& = \sum_{c=0}^{p-1} \sum_{a,b = 0}^M
    \hat\sigma(a - b)\ip{f_{c,b}}{f_{c,a}} \ip{B\be_{bp + c}}{B\be_{ap + c}} \\
& = \sum_{c=0}^{p-1} \int_{\TT}\,\sum_{a,b = 0}^M \gamma^{b-a}
    \ip{f_{c,b} \otimes B\be_{bp + c}}{f_{c,a} \otimes B\be_{ap + c}}
    \,d\sigma(\gamma) \\
& = \sum_{c=0}^{p-1} \int_{\TT}\,
    \biggl\langle\sum_{b=0}^M \gamma^b f_{c,b} \otimes B\be_{bp + c},
       \sum_{a=0}^M \gamma^a f_{c,a} \otimes B\be_{ap + c}\biggr\rangle
    \,d\sigma(\gamma) \\
& = \sum_{c=0}^{p-1} \int_{\TT}\,
    \biggl\lVert\sum_{a=0}^M \gamma^a f_{c,a}
          \otimes B\be_{ap + c}\biggr\rVert^2 \,d\sigma(\gamma) \\
& \ge 0, 
\end{align*}
so $\omega$ is positive.

To see that $\omega$ is decreasing, suppose $B\in\bs_m$.
Then $\beta(B) \in\bs_{m+1}$, so
\begin{align*}
\omega\circ\beta(B)
& = \sum_{k,l = 0}^{m+1} \lambda_{k,l}\ip{\beta(B)\be_l}{\be_k} \\
& = \sum_{k,l = 0}^{m+1} \lambda_{k,l}
    \sum_{i=1}^n \ip{r(v_i)B r(v_i)^*\be_l}{\be_k} \\
& = \sum_{k,l = 1}^{m+1} \lambda_{k,l}
    \sum_{i=1}^n \ip{e_l}{v_i}\ip{v_i}{e_k}\ip{B\be_{l-1}}{\be_{k-1}} \\
& = \sum_{k,l = 1}^{m+1} \lambda_{k,l}\ip{e_l}{e_k}\ip{B\be_{l-1}}{\be_{k-1}} \\
& = \sum_{k,l = 0}^m \lambda_{k+1,l+1}\ip{e_{l+1}}{e_{k+1}}\ip{B\be_l}{\be_k} \\
& = \omega(B)
\end{align*}
since from \eqref{eq:coeff} it is evident that
$\lambda_{k,l} = \lambda_{k+1,l+1}\ip{e_{l+1}}{e_{k+1}}$ for every $k,l$.

Let $\rho_\sigma = \Delta\omega\circ l$.
If $x\in\Ee_k$ and $y\in\Ee_l$, then
\[
\rho_\sigma(xy^*)
= \omega(\rankone xy)
= \lambda_{k,l}\ip{(\rankone xy)\be_l}{\be_k}
= \lambda_{k,l}\ip x{\be_k}\ip{\be_l}y,
\]
so $\rho_\sigma$ satisfies \eqref{eq:rhosigma2} as claimed.

Since the Fourier transform is linear, it is evident from the
defining formulas \eqref{eq:rhosigma} and \eqref{eq:coeff} that
the map $\sigma \mapsto \rho_\sigma$ is affine.
It remains only to show that $\sigma\mapsto\rho_\sigma$ is injective.
For this, suppose $a\in\ZZ$.  Since we have phased the sequence $(e_i)$
so that $\ip{e_i}{e_{i+p}}$ is always nonnegative,
by \cite[Proposition~1.2]{gui} the assumption that the series
$\sum\abs{1 - \ip{e_i}{e_{i+p}}}$ converges
implies that $\sum \abs{1 - \ip{e_i}{e_{i+ap}}}$ also converges.
Consequently the infinite product $\prod \ip{e_i}{e_{i+ap}}$ converges;
that is, there is a positive integer $i_a$ such that
$\prod_{i = i_a}^\infty \ip{e_i}{e_{i + ap}}$ exists and is nonzero.
Since
\[
\rho_\sigma(\be_{i_a + ap - 1}\be_{i_a - 1}^*)
 = \hat\sigma(a) \prod_{i = i_a}^\infty \ip{e_i}{e_{i + ap}}
\]
and the Fourier transform is injective, this shows that
$\sigma\mapsto\rho_\sigma$ is injective.
\end{proof}

\begin{cor}\label{cor:transitive}
Suppose $\rho$ is a pure essential product state of $\F$
with finite period $p$.  Then the gauge group acts $p$-to-$1$ transitively
on the extensions of $\rho$ to pure states of $\To$.
In particular, $\rho$ has precisely a circle of extensions to pure
states of $\To$.
\end{cor}

\begin{proof}
Suppose $\tilde\rho$ is an extension of $\rho$ to a pure state of $\To$.
Then there is a Borel probability measure $\sigma$ on $\TT$ such that
$\tilde\rho$ is the extension $\rho_\sigma$ of Theorem~\ref{theorem:extensions}.
Moreover, since $\tilde\rho$ is pure, there is a $z\in\TT$ such that
$\sigma$ is the point measure at $z$.
If $\lambda\in\TT$, then the pure state $\tilde\rho\circ\gamma_\lambda$
is equal to $\rho_\phi$, where $\phi$ is the point measure at $\lambda^p z$.
Thus the gauge group acts $p$-to-$1$ transitively
on the extensions of $\rho$ to pure states of $\To$.
\end{proof}

\begin{remark}
We conjecture that Corollary~\ref{cor:transitive} holds more generally for
non-product states.
\end{remark}


\begin{thebibliography}{20}

\bibitem{ace} H. Araki, A.L. Carey and D.E. Evans, {\em On $O_{n+1}$},
  J. Operator Theory {\bf 12} (1984), 247--264.

\bibitem{arv} W. Arveson, {\em Continuous analogues of Fock space},
Memoirs Amer. Math. Soc. {\bf 80} (1989), No. 409.

\bibitem{arv2} \bysame, {\em Continuous analogues of Fock space II:
the spectral \cstar algebra}, J. Funct. Anal. {\bf 90} (1990), no. 1, 138--205.

\bibitem{arv3} \bysame, {\em Continuous analogues of Fock space III:
singular states}, J. Operator Theory, J. Operator Theory {\bf 22} (1989),
no. 1, 165--205.

\bibitem{arv4} \bysame, {\em Continuous analogues of Fock space IV:
essential states}, Acta math. {\bf 164} (1990), 265--300.

\bibitem{bj} O. Bratteli and P.E.T. Jorgensen,
{\em Endomorphisms of $\Bb(\Hh)$, II: Finitely correlated states on $\Oo_n$},
J. Funct. Anal. (to appear).

\bibitem{bjp} O. Bratteli, P.E.T. Jorgensen and G. L. Price,
{\em Endomorphisms of $\Bb(\Hh)$}, Quantization, nonlinear partial differential
equations, and operator algebra (Cambridge, MA, 1994), 93--138,
Proc. Sympos. Pure Math., {\bf 59}, Amer. Math. Soc., Providence, RI, 1996

\bibitem{cun} J. Cuntz, {\em  Simple \cstar algebras generated by isometries},
  Comm. Math. Phys.   {\bf 57} (1977), 173--185. 

\bibitem{evans} D.E. Evans, {\em On $O_n$}, Publ. RIMS, Kyoto Univ.
{\bf 16} (1980), 915--927.

\bibitem{fowrae} N. Fowler and I. Raeburn, {\em Discrete product systems
and twisted crossed products by semigroups}, preprint, Univ. of Newcastle,
May 1997.

\bibitem{glimm} J. Glimm, {\em On a certain class of operator algebras},
Trans. Amer. Math. Soc. {\bf 95} (1960), 318--340.

\bibitem{gui} A. Guichardet, {\em Produits tensoriels infinis et
re\-pr\'e\-sen\-ta\-tions des relations d'anti\-commu\-ta\-tion},
Ann. Ec. Norm. Sup. {\bf 83} (1966), 1--52.

\bibitem{laca1} M. Laca, {\em Endomorphisms of $B(H)$ and Cuntz algebras},
  J. Operator Theory {\bf 30} (1993), 85--108.

\bibitem{laca2} \bysame, {\em Gauge invariant states of $\Oo_\infty$},
  J. Operator Theory {\bf 30} (1993), 381--396.

\bibitem{powers} R. T. Powers, {\em Representations of uniformly hyperfinite
algebras and their associated von Neumann rings}, Ann. Math. {\bf 86}
(1967), 138--171.

\bibitem{vonNeumann} J. von Neumann, {\em On infinite direct products},
  Comp. Math. {\bf 6} (1938), 1--77. 

\end{thebibliography}
\end{document}